\documentclass
[11pt]{article}

\usepackage{hyperref}
\usepackage{amsmath,amsfonts}

\usepackage{amsthm}
\usepackage{vmargin}
\usepackage{comment}
\includecomment{nonblinded}
\setpapersize{USletter}
\setmargrb{20mm}{20mm}{20mm}{20mm}
\title{Sequential Implementation of Monte Carlo Tests 
 with Uniformly Bounded Resampling Risk}
\author{Axel Gandy\\
\small Department of Mathematics,
\small Imperial College London, United Kingdom\\
a.gandy@imperial.ac.uk}
\usepackage{natbib}
\usepackage{graphics}

\newtheorem{theorem}{Theorem}
\newtheorem{lemma}[theorem]{Lemma}

\newcommand{\Prob}[1]{\text{\textup P}_{\!#1}}
\newcommand{\Expec}[1]{\text{\textup E}_{#1}}
\newcommand{\N}{\mathbb{N}}
\newcommand{\Z}{\mathbb{Z}}
\newcommand{\ind}[1]{1\!\!\!\:\textup{I}\left\{#1\right\}}
\DeclareMathOperator{\RR}{RR}
\DeclareMathOperator{\E}{E}

\begin{document}
\maketitle

\begin{abstract}

  \normalsize This paper introduces an open-ended sequential algorithm
  for computing the $p$-value of a test using Monte Carlo simulation.
  It guarantees that the resampling risk, the probability of a
  different decision than the one based on the theoretical $p$-value,
  is uniformly bounded by an arbitrarily small constant. Previously
  suggested sequential or non-sequential algorithms, using a bounded
  sample size, do not have this property. Although the algorithm is
  open-ended, the expected number of steps is finite, except when the
  $p$-value is on the threshold between rejecting and not rejecting.
  The algorithm is suitable as standard for implementing tests that
  require (re-)sampling.  It can also be used in other situations: to
  check whether a test is conservative, iteratively to implement
  double bootstrap tests, and to determine the sample size required
  for a certain power.

\end{abstract}


\section{Introduction}

Consider a statistical test that rejects the null hypothesis for large
values of a test statistic $T$.  Having observed a realization $t$,
one usually wants to compute the $p$-value 
$$
p=\Prob{}(T\geq t)
$$
where, ideally, $\Prob{}$ is the true probability measure under the
null hypothesis.  Of course, when the null hypothesis is composite,
$\Prob{}$ is often estimated (parametrically or non-parametrically).

In many cases, e.g.\ for bootstrap tests, the $p$-value $p$ cannot be
evaluated explicitly.  The usual remedy is a Monte Carlo test that
essentially replaces $p$ by
$$\hat p_{\text{naive}}=\frac{1}{n}\sum_{i=1}^n \ind{T_i\geq t},$$
where $T_1,\dots, T_n$ are independent replicates of the test
statistic $T$ under $\Prob{}$ and $\ind{}$ denotes the indicator function.

A reasonable requirement for a statistical method is what
\cite{gleser96:_commen_boost_confid_inter_t} called the \emph{first
  law of applied statistics}: ``Two individuals using the same
statistical method on the same data should arrive at the same
conclusion.''  For a test, ``conclusion'' is whether it rejects or
not, i.e.\ whether $p$ is above or below a given threshold $\alpha$
(often $\alpha=0.05$).

Monte Carlo tests do not satisfy this law: For an estimator $\hat p$,
let the {\it resampling risk} $\RR_p(\hat p)$ be the probability that
$\hat p$ and the true $p$ are on different sides of the threshold
$\alpha$.  More precisely,
 $$\RR_p(\hat p)\equiv
\begin{cases}
   \Prob{p}(\hat p> \alpha)&\text{if }p\leq \alpha,\\
   \Prob{p}(\hat p\leq \alpha)&\text{if }p> \alpha.
 \end{cases}$$ 
The resampling risk $\RR_p(\hat p_{naive})$ of the
 naive estimator $\hat p_{\text{naive}}$ can be substantial, e.g.\
 $\RR_p(\hat p_{naive})=0.146$ for $n=999$, $\alpha=0.1$, and $p=0.11$
 \cite[Table 1]{joeckel84:_comput_aspec_of_monte_carlo_tests}.
 Furthermore, no matter how large $n$ is chosen,
 $$\sup_{p\in [0,1]} \RR_p(\hat p_{naive})\geq 0.5.$$

 In the present article, we introduce a recursively defined sequential algorithm which
 gives an estimator $\hat p$ of $p$ that uniformly bounds the
 resampling risk as follows:
$$
\sup_{p\in [0,1]} \RR_p(\hat p)\leq \epsilon
$$
for some arbitrary (small) $\epsilon>0$.
Although the algorithm is open-ended, i.e.\ the number of steps is not
bounded, the expected number of steps is finite for $p\neq \alpha$.
In particular, if $p$ is far away from the threshold $\alpha$ then the
algorithm usually stops quickly.

Having reached step $n$ without stopping, there exists an interval (with
length going to $0$ as $n\to \infty$) that contains the not yet available
estimate $\hat p$. This interval can be used as an interim result.
  
It is well known, that Monte Carlo tests may lose  power compared
to the theoretical test, see e.g.\ \cite{Hope1968SMC} or \citet[p.\
155,156]{davison1997bma}.  Our algorithm bounds this loss of
power by the arbitrarily small constant $\epsilon$.

The proposed sequential algorithm can be used as standard
implementation for (re-)sampling based tests in statistical software.
Essentially, one only has to set $\epsilon$ to a suitably small default
value (e.g.\ $10^{-3}$ or $10^{-5}$), and ensure that
the algorithm reports intermediate results until it finishes.
An \emph{R} package is available from the author's web page

\begin{nonblinded}
\ \url{http://www.ma.ic.ac.uk/~agandy}
\end{nonblinded}
.

Other sequential procedures to compute $p$-values have been suggested
previously.  The suggestion of \cite{Davi:MacK:boot:2000} is
relatively close to our algorithm.  They also use a uniform bound on
the resampling risk as motivation.  However, their algorithm does not
really guarantee this bound since, when deciding whether to stop,
 they do not take into account the problem
of multiple testing.  Furthermore, whereas we allow stopping after
each step, they only allow stopping after $2^kB$ steps for
$k=0,\dots,n$, where $B$ is some constant.

\cite{besagclifford91:_sequen_monte_carlo} suggest a sequential
procedure which stops if the partial sum $\sum_{i=1}^n \ind{T_i\geq t}$ reaches a
given threshold or if a given number of samples $n$ is reached. The
motivation is that if $p$ is high, i.e.\ if the test result is far away
from being significant, fewer replicates are needed than in the naive
approach.  More recently, \cite{Fay2007UTS} suggested using a truncated
sequential probability ratio test.

\cite{Andr:Buch:thre:2000,Andr:Buch:eval:2001} suggest using as
criterion the relative difference between the $p$-value from the
finite-sample bootstrap and the 'ideal' bootstrap using infinite
sample size.  This method involves drawing some fixed number of
bootstrap samples and then using asymptotic arguments to determine the
number of repetitions needed.  Once the number of repetitions has been
chosen, the test can be performed by drawing the remaining bootstrap
repetitions (without further sequential considerations).

Bayesian approaches have been suggested in \cite{Lai:near:1988} and in
\cite{fay02:desig_monte_carlo}.  By putting a (prior) distribution on
$p$, the average resampling risk $\E(\RR_p(\hat p))$ can be bounded.
This bound is much weaker than our uniform bound on $\RR_p$.

The present article is structured as follows: In Section
\ref{sec:algorithm}, we precisely define the sequential algorithm and
describe the key results of the paper.  In Section
\ref{sec:some-remarks}, we comment on several aspects of our algorithm
such as the expected number of steps, choice of tuning parameters, and
details of the implementation.

In Section \ref{sec:an-application}, we demonstrate the wide
applicability of our sequential algorithm in a simple practical
example.  Proofs are relegated to the appendix.

\section{The Algorithm and Key Results}
\label{sec:algorithm}

Instead of considering independent replicates of $\ind{T\geq t}$ one
can obviously consider replicates from a Bernoulli distribution with
parameter $p$.  From now on, let $X_1, X_2, \dots $ be independent and
identically distributed Bernoulli distributed random variables with
parameter $p$.  In the notation of the introduction, $X_i=\ind{T_i\geq
  t}$.  Expectations and probabilities are taken using the $p$ that is
indicated in a subscript ($\Prob{p}(\cdot)$ and $\Expec{p}(\cdot)$).

\begin{figure}[tbp]
  \centering
\includegraphics{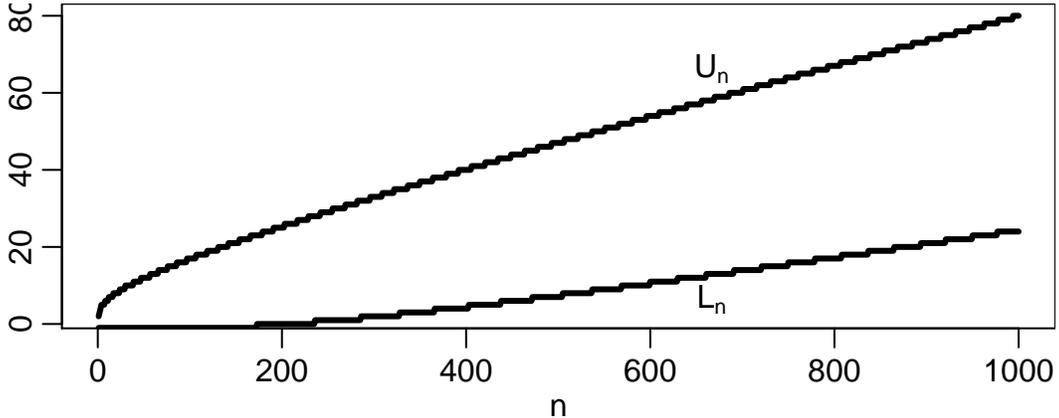}
  \caption{The algorithm stops if $S_n\geq U_n$ or
$S_n\leq L_n$
The boundaries are computed for the threshold 
$\alpha=0.05$ using the spending sequence
  $\epsilon_n=\frac{n}{1000(1000+n)}$.}
  \label{tab:boundaries}
\end{figure}

Our sequential algorithm stops once the partial sum  $S_n = \sum_{i=1}^n
X_i$ hits
boundaries given by two integer sequences $(U_n)_{n\in \N}$ and
$(L_n)_{n\in \N}$ with $U_n\geq L_n$, i.e.\ we stop after
$$
\tau = \inf \{k\in \N: S_k\geq U_k\text{ or } S_k \leq L_k\}
$$
steps. In the above, $\N=\{1,2,\dots\}$.  
 Figure \ref{tab:boundaries} shows an example of 
sequences $U_n$ and $L_n$ resulting from the following definition.

We construct $U_n$ and $L_n$ such that for $p\leq \alpha$ (resp.\ $p>
\alpha$) the probability of hitting the upper boundary $U_n$ (resp.\
lower boundary $L_n$) is at most $\epsilon$, where $\epsilon>0$ is the
desired bound on the resampling risk. We will see that it suffices to
ensure that for $p=\alpha$ the probability of hitting each boundary is
at most $\epsilon$.  We use a recursive definition that for each $n$
minimizes $U_n$ (resp.\ maximizes $L_n$) conditional on
\begin{align}
\label{eq:hitputuntilk}
\Prob{\alpha}(\text{hit upper boundary until }n)=\Prob{\alpha}(\tau \leq n, S_\tau\geq U_\tau)&\leq \epsilon_n \text{ and}\\
\label{eq:hitdownuntilk}
\Prob{\alpha}(\text{hit lower boundary until }n)=\Prob{\alpha}(\tau \leq n, S_\tau\leq L_\tau)&\leq \epsilon_n,
\end{align}
where $\epsilon_n$ is a non-decreasing sequence with $\epsilon_n \to
\epsilon$ and $0\leq \epsilon_n < \epsilon$. We call 
$\epsilon_n$ \emph{spending sequence}.
The sequence $\epsilon_n$ is used to control how fast the allowed
resampling risk $\epsilon$ is spent. 
In the examples of this paper, we will
use $\epsilon_n=\epsilon\frac{n}{k+n}$ for some constant $k$.  
There is a close connection of our spending function $\epsilon_n$ to the
$\alpha$-spending function (or ``use'' function) of
\cite{Lan:DeMe:disc:1983}. 

The formal definition of  the boundaries  $(U_n)$ and $(L_n)$  is as 
follows: Let  $U_1 = 2, L_1 = -1$ and, recursively for $n\in \{2,3,\dots\}$, let
\begin{equation}
\begin{split}
\label{defUnLnsimple}
U_n=&\min\{j\in \N: \Prob{\alpha}(\tau \geq n,S_n\geq j)+ \Prob{\alpha}( \tau <n, S_{\tau}\geq U_\tau)\leq \epsilon_n\},  \\
L_n=&\max\{j\in \Z:   \Prob{\alpha}(\tau \geq n,S_n\leq j )+\Prob{\alpha}(\tau <n, S_\tau \leq L_\tau)\leq \epsilon_n
 \}.
\end{split}
\end{equation}
Note that $U_n$ (resp.\ $L_n$) is the minimal (resp.\ maximal) 
value for which (\ref{eq:hitputuntilk}) (resp.\ (\ref{eq:hitdownuntilk}))
holds true given $U_1,\dots, U_{n-1}$ and $L_1,\dots, L_{n-1}$.
Using induction, one can 
see that (\ref{eq:hitputuntilk}) and
(\ref{eq:hitdownuntilk}) hold true for all $n$.

The following theorem shows that 
the expected number of steps of the algorithm is finite for $p\neq \alpha$ 
and that the probability of hitting the
``wrong'' boundary is bounded by $\epsilon$.
\begin{theorem}
\label{th:stopinfinitetime}
 Suppose that $\epsilon<1/2$ and  $\log(\epsilon_n - \epsilon_{n-1})=o(n)$ as $n\to \infty$. 
Then $U_n-\alpha n=o(n)$, $\alpha n - L_n=o(n)$ 
and $\Expec{p}(\tau)<\infty$ for all $p\neq \alpha$.
Furthermore,
\begin{equation}
  \label{eq:boundhitwrongboundary}
  \sup_{p\in [0,\alpha]}\Prob{p}(\tau<\infty,S_\tau \geq U_\tau)\leq \epsilon
\text{ and } \sup_{p\in (\alpha,1]}\Prob{p}(\tau<\infty,S_\tau \leq L_\tau)\leq \epsilon.
\end{equation}

\end{theorem}
\noindent The proof of this and the next
 theorem can be found in the appendix.

\begin{figure}[tbp]
  \centering
\includegraphics{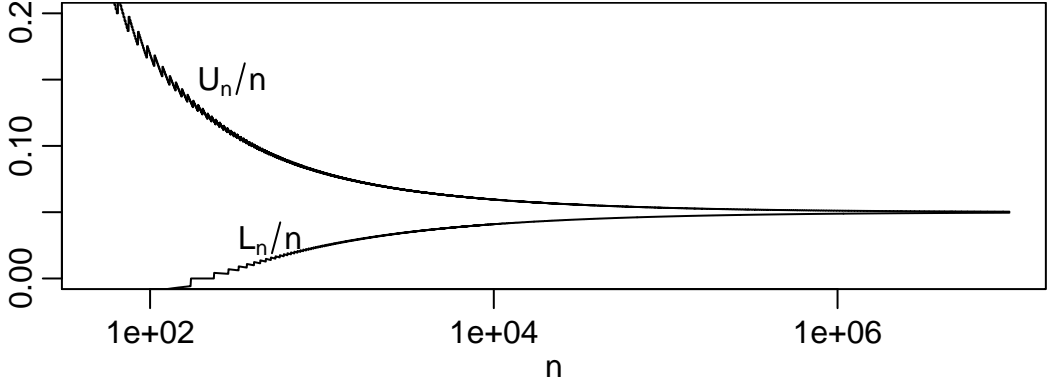}
  \caption{If the upper boundary is hit after
 $\tau$ steps then $\hat p=U_\tau/\tau$, if the lower boundary is hit then $\hat p=L_\tau/\tau$.
   We used $\alpha=0.05$ and $\epsilon_n= \frac{n}{1000(1000+n)}$. Note the log-scale on the horizontal axis.
  }
  \label{fig:monotonicity}
\end{figure}

As estimator for $p$ we  use the maximum likelihood estimator
$$
\hat p = 
\begin{cases}
\displaystyle\frac{S_\tau}{\tau},&\tau<\infty,\\
\alpha, & \tau = \infty.  
\end{cases}
$$
Figure \ref{fig:monotonicity} shows how the estimator depends on
when the boundary is hit.

The next theorem gives the uniform bound on the resampling risk.
\begin{theorem}
\label{th:MLE}
Suppose $\epsilon\leq 1/4$ and  $\log(\epsilon_n - \epsilon_{n-1})=o(n)$ as $n\to \infty$. Then
$$
\sup_{p\in[0,1]} \RR_p(\hat p)\leq \epsilon.
$$
\end{theorem}
Note that our default spending sequence $\epsilon_n=\epsilon\frac{n}{k+n}$ satisfies
the conditions of the above theorems.

The conditions in the above two theorems are not minimal.  For example,
considering whether to stop  only every $\nu\in \N$ steps will,
using a slightly modified proof, lead to the same results.

\section{Remarks}
\label{sec:some-remarks}

\subsection{Lower Bound on the Expected Number of Steps}
\label{sec:lowerboundonsteps}
Suppose $\hat p$ is a sequential algorithm with stopping time
$\tau$ that has a uniformly bounded resampling risk
$$
\sup_{p\in [0,1]}RR_p(\hat p)\leq \epsilon.
$$
We derive a lower bound on $\Expec{p}(\tau)$ and show that 
in a Bayesian setup the expected number of steps 
is infinite.

For $p_0>\alpha$ consider the hypotheses 
$H_0: p=\alpha$ against $H_1:p=p_0$.
We can construct a test by rejecting $H_0$ iff $\hat p>\alpha $.
The probability of both the type I and the type II error is $\epsilon$.

Consider the sequential probability ratio test (SPRT), see \cite{wald1945STo},
 of $H_0$ against
$H_1$ with the same error probabilities. Let $\sigma_{p_0}$ denote its
stopping time. The SPRT minimizes the expected number of steps among
all sequential tests with the same error probabilities
\citep{Wald1948OCo}.  Thus, 
\begin{align}
\label{lowerboundrunningtime}
\Expec{p_0}(\tau)\geq \Expec{p_0}(\sigma_{p_0})
\approx\frac{(1-\epsilon)\log\left(\frac{1-\epsilon}{\epsilon}\right)
+\epsilon\log\left(\frac{\epsilon}{1-\epsilon}\right)}
{p_0\log(\frac{p_0}{\alpha})+(1-p_0)\log\left(\frac{1-p_0}{1-\alpha}\right)},
\end{align}
where the approximation is from \cite[(4.8)]{wald1945STo}.  The
approximation can be replaced by an inequality via
\citet[(4.13),(4.15) or (4.16)]{wald1945STo}.

Equation (\ref{lowerboundrunningtime}) also holds true for $p_0<\alpha$.
Indeed, one only needs to replace the above $H_0$ by
$H_0:p=\alpha+\delta$ for some $\delta>0$ and let $\delta\to 0$.

Suppose, in a Bayesian sense, that $p$ is random,  
having distribution function $F$ with derivative $F'(\alpha)>0$. 
Then for some $c>0$ and some $\delta>0$, 
\begin{align*}
\Expec{}(\tau)=\int_0^1 \Expec{p}(\tau)dF(p)\geq
c \int_{\alpha-\delta}^{\alpha+\delta} 
\left(p\log\left(\frac{p}{\alpha}\right)+(1-p)\log\left(\frac{1-p}{1-\alpha}\right)\right)^{-1}\!\!\!dp=\infty
\end{align*}
The last equality holds since the integrand is proportional to
$2\alpha(1-\alpha)(p-\alpha)^{-2}$ as $p\to \alpha$ (by  e.g.\ l'Hospital's
rule),.

Suppose we want to use our algorithm to compute the power or the level
of a test.  Then $p$ is indeed random and $\E (\tau)= \infty$.  Thus for
this application we have to truncate our stopping time (e.g.\ by some
deterministic constant).

\subsection{Error Bound and Spending Sequence}
\begin{figure}[tbp]
  \centering
\includegraphics{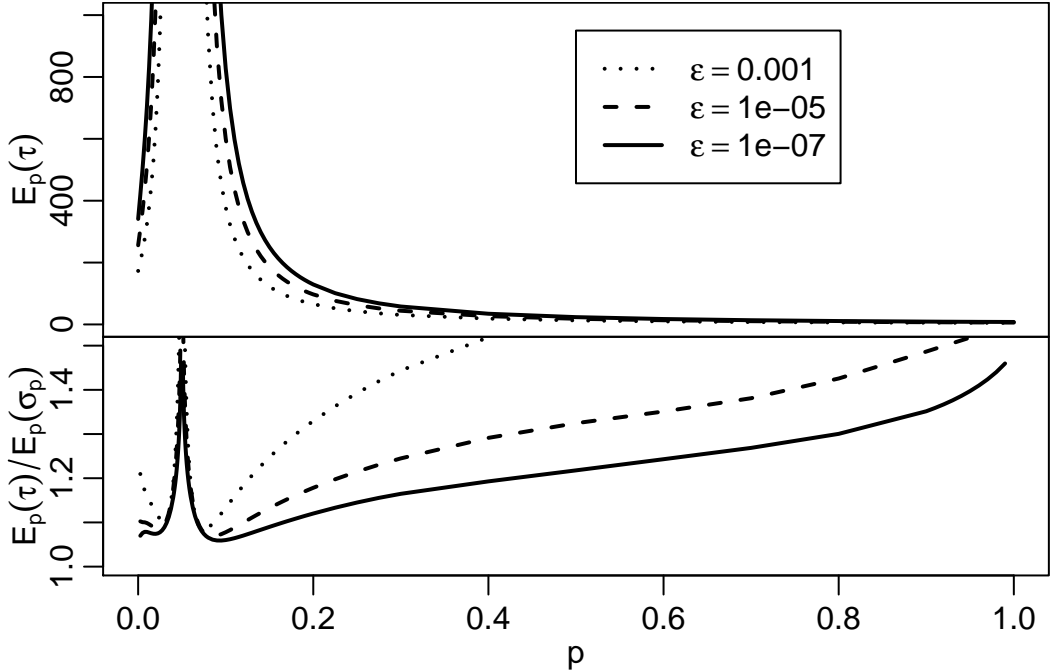}
  \caption{Top: Expected number of steps $\Expec{p}(\tau)$ of the algorithm against  the true parameter $p$.
Bottom:  $\Expec{p}(\tau)$ divided by the theoretical lower limit in (\ref{lowerboundrunningtime}).
    The threshold $\alpha=0.05$ and the
spending sequence  $\epsilon_n=\epsilon \frac{n}{1000+n}$ were used.
  }
  \label{tab:runningtime}
\end{figure}

Figure \ref{tab:runningtime} illustrates the dependence of the
expected number of steps $\Expec{p}(\tau)$ on the true $p$ and on the
error bound $\epsilon$.  For most $p$, the algorithm stops quite
quickly.  Furthermore, the dependence on the bound $\epsilon$ of the
resampling risk is only slight, so $\epsilon$ can be chosen small in

The lower part of Figure \ref{tab:runningtime} shows that our
algorithm with the default spending sequence is not too far away from
the theoretical boundary.  Can this be improved by choosing a
different $\epsilon_n$?  What should ``improved'' mean? There is no
obvious optimality criterion.  As Section \ref{sec:lowerboundonsteps}
shows, a criterion like the average number of steps under the null
hypothesis cannot be used since it is always infinite.  An option is
to try to minimize something like
$\int_0^1\Expec{p}(\tau)/\Expec{p}(\sigma_p)dp,$ i.e.\ integrating the
function plotted in the lower part of Figure \ref{tab:runningtime}.
However, pursuing this  further is beyond the scope of this
article.

In a small study, not reported here, we have looked at other choices
besides our default $\epsilon_n=\epsilon\frac{n}{k+n}$. The main
conclusion is that the choice of $\epsilon_n$ does not seem to have a
big influence - as long as the allowed error is spent at a
sub-exponential rate (satisfying the conditions in Theorems
\ref{th:stopinfinitetime} and \ref{th:MLE}).

\subsection{Bounds on the Estimator Before Stopping}
In practice, the algorithm should report back after a fixed number of
steps, even if it has not stopped yet.  After this intermediate stop
one can continue the algorithm.

In the case of such an intermediate stop one can compute an interval
in which $\hat p$ will eventually lie.  One can base this interval on
the inequalities (\ref{eq:asymptboundUn}) and (\ref{eq:asymptboundLn})
in the appendix. Indeed, after $n$ steps, conditional on $\tau\geq n$,
one gets $\hat p\in
[\alpha-\frac{\Delta_n+1}{n},\alpha+\frac{\Delta_n+1}{n}]$ for $\tau
\geq n$, where $\Delta_n= \sqrt{-n\log(\epsilon_n-\epsilon_{n-1})/2}$.
Note that under the assumptions of Theorem \ref{th:stopinfinitetime}
we have $(\Delta_n+1)/n\to 0$. 

 These bounds are not very tight. They can be improved as follows.
Conditional on not having stopped after $n$ steps, we have 
\begin{equation}
  \label{eq:boundsonestlookahead}
  \min_{\nu\geq n}\frac{L_\nu}{\nu} \leq \hat p \leq \max_{\nu\geq n} \frac{U_\nu}{\nu}
\end{equation}
(in the proof of Theorem \ref{th:MLE} we show $U_\nu\geq\nu\alpha\geq
L_\nu$).  Figure \ref{fig:monotonicity} shows a plot of
$\frac{1}{n}U_n$ and $\frac{1}{n}L_n$ for one particular spending
sequence.  It seems that $\frac{1}{\nu}U_\nu$ is overall decreasing,
and that $\frac{1}{\nu}L_\nu$ is overall increasing.  This also seems
to be true for further spending sequences of the type
$\epsilon_n=\epsilon\frac{n}{k+n}$.  Thus an ad-hoc way of computing
the upper bound in (\ref{eq:boundsonestlookahead}) is by evaluating $
\max_{n\leq \nu\leq n+\mu} \frac{U_\nu}{\nu}$ where $\mu$ is chosen
suitably, e.g.\ $\mu=2/\alpha$. A similar argument applies, of course,
to the lower bound of (\ref{eq:boundsonestlookahead}).

Instead of reporting back after a fixed number of steps, the algorithm
could also report back after a certain computation time, e.g.\ one
minute.  This has the advantage that a sensible default value can be
used irrespective of the time one sampling step takes.

\subsection{Confidence Intervals}

Confidence intervals for $p$ can be constructed similarly to 
\cite{Armitage1958NSi}.  Suppose the
algorithm stops and returns $p^{obs}$ as result.  Then a $1-\beta$ 
confidence interval for $p$ is given by $[\underline p,\overline p]$, where
$\underline p=0$ for $p^{obs}=0$, $\overline p=1$ for $p^{obs}=1$, and otherwise
\begin{align*}
   \Prob{\underline p}(\hat p\geq p^{obs})= \beta/2,\quad
   \Prob{\overline p}(\hat p\leq p^{obs})= 1- \beta/2.
\end{align*}
\cite{Armitage1958NSi} showed that the probabilities on the left
hand sides are strictly monotonic in $\underline p$ and $\overline p$ and thus 
 $\underline p$ and $\overline p$ are well-defined.

One can compute $\underline p$ and $\overline p$ numerically.
If the computation of the above probabilities involves an infinite sum
we consider the complement event instead, e.g.\ we may replace $\Prob{\underline p}(\hat p\geq p^{obs})$ 
by  $1-\Prob{\underline p}(\hat p< p^{obs})$.

Suppose the algorithm has not stopped, i.e.\ $\tau > n$.  Then by the
arguments of the previous subsection we get $\hat p_{\min}$, $\hat
p_{\max}$ such that $\hat p\in [\hat p_{\min}, \hat p_{\max}]$.
Replacing $p^{obs}$ in the definition of $\underline p$ (resp.\ $\overline p$) by $\hat
p_{\min}$ (resp.\ $\hat p_{\max}$) 
produces an interval that includes the
confidence interval one gets once the algorithm has finished.  Thus
this is a confidence interval itself, with a (slightly) increased coverage
probability.

\subsection{Implementation Details}

To compute $U_n$ and $L_n$  via (\ref{defUnLnsimple}), 
one needs to know the distribution of $S_n$ given $\tau \geq n$ as
well as $\Prob{\alpha}(\tau< n, S_\tau \geq U_\tau)$ and
$\Prob{\alpha}(\tau< n, S_\tau\leq L_\tau)$.  These quantities can be
updated recursively.  Furthermore, the amount of memory required to
store these quantities is proportional to $U_n - L_n$.

\begin{figure}[tbhp]
\centering
\includegraphics{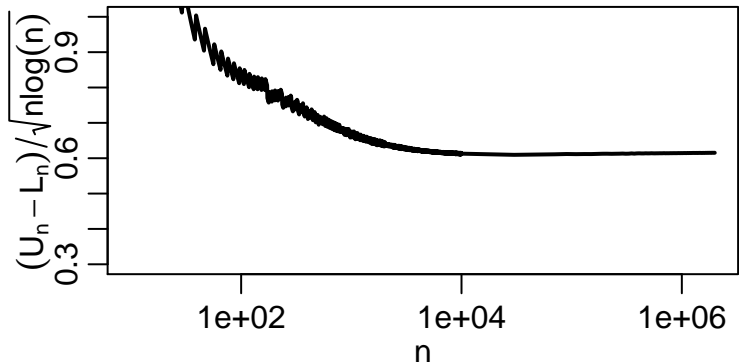}
\caption{ 
$U_n-L_n$ seems to be  roughly proportional to $\sqrt{n  \log n}$.
We used  $\alpha =0.05$ and
 $\epsilon_n  = \frac{n}{1000(1000+n)}$.
Note the log-scale on the horizontal axis.
}
\label{fig:boundariesasympt}
\end{figure}

What is the additional computational effort for the sequential
procedure?  The main effort at each step is to compute the
distribution of $S_n$ given $\tau \geq n$ from the distribution of
$S_{n-1}$ given $\tau \geq n-1$. This effort is proportional to $U_n
-L_n$.  Hence, if the sequential procedure stops after $n$ steps the
computational effort is roughly proportional to
$\sum_{i=1}^n|U_i-L_i|$.  In order to get an idea of how big $U_n -
L_n$ is we considered a specific example in Figure
\ref{fig:boundariesasympt}.  In this example it seems as if $U_n - L_n
\sim \sqrt{n\log n}$.  The overhead of the algorithm can  be removed
through precomputation of $L_n$ and $U_n$.

Our sequential procedure can be easily parallelized, 
e.g.\ by distributing the generation  of the samples.

\subsection{Using the Algorithm as a Building Block}
\label{sec:iterat}

Our procedure can be used as a building block in more complicated
computations. For example, one can use the sequential procedures in
this paper to estimate the power of a resampling based test, by using
the algorithm in the ``inner'' loop.  Because of the problems
mentioned in Section \ref{sec:lowerboundonsteps}, the number of
replications in the inner loop have to be restricted by a constant.
Of course, this is rather  ad-hoc, but it should give a
similar performance (with less computational effort) than the naive
approach of nesting two loops within one another.

For the problem of computing the power of a bootstrap test, some
dedicated algorithms exists, such as that suggested by
\cite{Boos:Zhan:mont:2000}.  Their algorithm can be combined with ours
by using our sequential procedure in the inner loop.

To compute the power of a bootstrap test,
\cite{jennison92:_boots_tests_and_confid_inter} has suggested a
sequential procedure for the ``inner'' loop.
\cite{jennison92:_boots_tests_and_confid_inter} uses an approximation
to bound the probability of deciding differently than the bootstrap
that uses only a fixed number of samples.  In contrast to that, the
present article bounds the probability of deciding differently than
the ``ideal'' bootstrap based on an infinite sample size.

Furthermore, the algorithm can be used iteratively, e.g.\ for double
bootstrap tests.  Examples  can be found in
Section \ref{sec:an-application}.

\section{Applications}
\label{sec:an-application}

This section demonstrates the wide applicability  of our algorithm
in a simple example,
  already used by \cite{mehta1983nap},
by \cite{newton1994brm}, and by
\citet[Example 4.22]{davison1997bma}.
Suppose  $39$ observations have been categorized
according to two categorical variables resulting
in counts given by the following 
 two-way sparse contingency table:
\begin{center}
\begin{tabular}{ccccccc}
  1&2&2&1&1&0&1\\
  2&0&0&2&3&0&0\\
  0&1&1&1&2&7&3\\ 
  1&1&2&0&0&0&1\\ 
  0&1&1&1&1&0&0
\end{tabular}
\end{center}
Let $A=(a_{ij})$ denote this  matrix.

Consider the test of the null hypothesis that  the two variables are
independent which rejects for large values of the likelihood ratio
test statistic
$$T(A)=2\sum_{i,j} a_{ij} \log(a_{ij}/h_{ij}),$$
where $h_{ij}=\sum_\nu a_{\nu j}\sum_\mu a_{i\mu}/\sum_{\nu\mu} a_{\nu
  \mu}$.  It is well known that under the null hypothesis, as the
sample sizes increases, the distribution of $T(A)$ converges to a
$\chi^2$-distribution with $(7-1)(5-1)=24$ degrees of freedom.
Applying this test to the above matrix leads to a $p$-value of
$0.031$.

\subsection{Parametric Bootstrap}
\label{sec:parametric-bootstrap}
Since the contingency table $A$  is sparse, the asymptotic approximation may
be  poor.
To remedy this, \cite{davison1997bma} suggested  a parametric bootstrap
that simulates under the null hypothesis  based on the 
row and column sums of $A$.

Using the naive test statistic $\hat p_{naive}$ with   $n=1{,}000$ replicates
results in a $p$-value of 
$0.041$. This is  below the usual threshold of 5\
thus the test would be interpreted as significant.
However, as further computations show, the probability of 
reporting a $p$-value larger than 5\

Next, we applied our algorithm using $\alpha=0.05$, $\epsilon=10^{-3}$
and $\epsilon_n=\epsilon\frac{n}{1000+n}$.  We shall use this
$\epsilon$ and $\epsilon_n$ in all other examples of Section
\ref{sec:an-application}.  Assume that we decide to let our algorithm
run for at most 1{,}000 steps initially.  Not having reached a
decision, the algorithm tells us that the final estimate will be in
the interval $[\hat p_{\min},\hat p_{\max}]=[
0.027,0.080
]$.  Our algorithm finally stops after $8{,}574$
samples, reporting a $p$-value of $0.040$.
The advantage of our algorithm is that we can be (almost)
certain that the ideal bootstrap would also return a significant
result.

\subsection{Some Notation}

To describe further uses of our algorithm we introduce the following
notation.  Let $h_\alpha$ be the function that applies our algorithm
with the threshold $\alpha$ to a sequence with elements in $\{0,1\}$
and returns the resulting estimate $\hat p$.  If the sequence is
finite, say of length $n$, and the algorithm has not stopped after $n$
steps then $h_\alpha$ simply returns the current estimate
$S_n/n$.

With this, the above use of our algorithm for the parametric bootstrap
can be written as
 $$h_{0.05}\left(\ind{T(A_i)\geq T(A)}_{i\in \N}\right),$$
where  $A_1, A_2, \dots$ denote independent samples
under the null hypothesis estimated from the row and column sums of
the matrix $A$.

\subsection{Checking the Level of Tests}
\label{sec:checking-level-tests}

As mentioned earlier, the asymptotic $\chi^2$-distribution may not be
a good approximation because the observed matrix $A$ is relatively
sparse.  We can use our algorithm to check whether the test is
conservative or liberal. To do this  at the  5\
 we estimate the
rejection rate by
 $$h_\alpha\left(\ind{T(A_i)\geq \chi^2_{24,0.95}}_{i\in \N}\right),$$
where  $\chi^2_{24,0.95}$ denotes the 0.95 quantile of the $\chi^2$-distribution
with 24 degrees of freedom.
We start our procedure with a threshold of $\alpha=0.05$.
It stops after $1{,}437$ steps and reports a $p$-value of 
$0.074$.
Hence, the test based on the asymptotic distribution  seems to be liberal.

How liberal is it?  To find out whether the rejection rate is above
$0.07$ we start our procedure with a threshold of $\alpha=0.07$.
After $66{,}736$ samples, the estimated rejection rate
is $0.075$. Thus it is (almost) certain that the
test at the nominal level 5\

Next, we check whether the parametric bootstrap test of Section
\ref{sec:parametric-bootstrap} does any better.  For this we use the
sequential procedure iteratively and compute the rejection rate by
$$
h_\alpha\left(\ind{h_{0.05}\left(\ind{T(A_{ij}))\geq
        T(A_i)}_{j=1,\dots,M}\right)\leq 0.05}_{i\in \N} \right),
$$
where for each $i\in \N$, $A_{i1}, A_{i2},
\dots$ denote independent samples under the null hypothesis estimated
from the matrix $A_i$.

As explained at the end of Section \ref{sec:lowerboundonsteps},
 in the ``inner'' use of $h_{0.05}$ we need to stop after a finite number $M$ of steps. For the following we use $M=250$.
Setting 
$\alpha=0.05$, the outer algorithm stops after $264$ steps yielding a $p$-value of 
 $0.114$.
Using $\alpha =0.07$ after  $1{,}769$ steps we get a $p$-value of $0.096$.
Hence, the bootstrap test seems to be quite liberal as well.

For $\alpha=0.05$ (resp.\ $\alpha=0.07$) we generated a total of
$21{,}250$ (resp.\ $131{,}552$) samples
in the inner loop.  A naive alternative consists of just two nested
loops.  To get a similar precision one could use $M$ steps in the
inner loop and $1{,}000$ steps in the outer loop.  For this,
$250{,}000$ samples need to be generated, far more than in our nested
sequential algorithm.

\subsection{Double Bootstrap}
\citet[Example 4.22]{davison1997bma} suggest that the parametric
bootstrap could be improved by using a double bootstrap.  The double
bootstrap employs two loops that are nested within one another.
\cite{davison1997bma} suggest that a sensible choice would be to
use roughly 1{,}000 steps in the outer loop and 250 steps in the inner
loop.  As the classical double bootstrap needs to resample once before
starting the inner loop, it needs $251{,}000$ resampling steps.

To reduce the number of steps,  we can use our algorithm iteratively:
First, compute the $p$-value from the parametric bootstrap using, say, 10{,}000 samples by
$$p=\frac{1}{10{,}000}\sum_{i=1}^{10{,}000} \ind{T(A_i)\geq T(A)}.$$
After that, we  compute
the $p$-value of the double bootstrap by
$$h_{0.05}\left(\ind{
       h_p\left(\ind{T(A_{ij})\geq T(A_i)}_{j=1,\dots,M}\right)
    \leq p}_{i\in \N}\right).$$

Applying this with $M=250$, the outer algorithm stops after $1{,}117$ samples
in the outer loop and returns a $p$-value of $0.078$,
which, in contrast to the previous tests, is not significant at the 5\
In the inner loop we used  only
$77{,}405$ samples.
Adding the $10{,}000$ samples needed to compute $p$ and
the $1{,}117$ samples from the null model fitted to $A$ the
total  number of samples generated is  $88{,}522$.
This  compares favorably to the $251{,}000$ for the classical double bootstrap.

To check the level of the double bootstrap test
we can combine the approach of Section \ref{sec:checking-level-tests}  with the approach 
of the current subsection. This results in iterating the procedure three-times.
For the double bootstrap we set $M=250$ and stop the outer algorithm
after  500 steps. 
If we check whether the true level of our algorithm at the asymptotic level
5\
and reports a $p$-value of $0.050$.
 Hence the double bootstrap seems to
be less liberal (if it is liberal at all) than the asymptotic test
or the simple parametric bootstrap.
In the innermost application of our algorithm we needed  $12{,}688{,}117$
 resampling steps.
The naive approach with a similar maximal number of steps for the inner loops and
1{,}000 steps for the outer loop  would have 
used $1{,}000\cdot 500\cdot 250 =1.25\cdot 10^8$ steps in the innermost loop, 
more than $9$ times  the number of samples our iterated
algorithm needed.

\subsection{Determining Sample Size}

In Section \ref{sec:checking-level-tests}, we have seen how to use our
algorithm to check the level of a test.  Similarly, with the obvious
modification of generating $A_i$ from the given alternative, one can
check whether a test achieves a desired power for a given sample size.
Furthermore, the the minimal sample size that achieves a certain power
can be found by combining our algorithm with e.g.\ a bisectioning
algorithm.

\section{Conclusions}

We presented a sequential procedure to compute $p$-values by sampling.
When the algorithm stops one has the ``peace of mind'' that, up to a
small error probability, the $p$-value reported by the procedure is on
the same side of some threshold as the theoretical $p$-value.  In
other words, the resampling risk is uniformly bounded by a small
constant. If the algorithm has not stopped then one can give an
interval in which the final estimate will be.

The basic algorithm can also be used in several other
situations.  It can be used to check whether a test is conservative or
liberal, in can be used iteratively for double bootstrap test, and it
can be used to determine the sample size needed to achieve a certain
power.

\begin{nonblinded} 
 \section*{Acknowledgment}

 Part of the work was carried out during a stay of the author at the
 Centre of Advanced Study at the Norwegian Academy
 of Science and Letters in Oslo.
\end{nonblinded}

\appendix

\section{Proofs}

The following lemma is needed in the proof of
 Theorem \ref{th:stopinfinitetime}.
\begin{lemma}
\label{le:monotonichitting}
For $0\leq p  \leq q\leq 1$,
\begin{align*}
\Prob{p}(\tau<\infty,S_\tau \geq U_\tau) \leq 
\Prob{q}(\tau<\infty,S_\tau \geq U_\tau),\\
\Prob{p}(\tau<\infty,S_\tau\leq L_\tau) \geq 
\Prob{q}(\tau<\infty,S_\tau\leq L_\tau).
\end{align*}
\end{lemma}
\begin{proof}
  Let  $V_1, V_2, \ldots $ be independent random variables with
a uniform distribution on $[0,1]$ under the probability measure $\Prob{\ast}$.
For $x\in [0,1]$, let
 $$S_{x,n}=\sum_{i=1}^n \ind{V_i\leq x}.$$ 
Clearly $S_{p,n}\leq S_{q,n}$.
Let $\tau_x=\inf\{k\in \N: S_{x,k}\geq U_k\text{ or } S_{x,k}\leq L_k\}$.
Then
\begin{equation}
\label{eq:ineqinlemma}
\begin{split}
\Prob{p}(\tau\leq n, S_\tau\geq U_\tau)
&= \Prob{\ast}(\tau_p\leq n, S_{p,\tau_p}\geq U_{\tau_p})\\
  &\stackrel{(\dag)}{\leq} \Prob{\ast}(\tau_q\leq n, S_{q,\tau_q}\geq U_{\tau_q})
=\Prob{q}(\tau\leq n, S_\tau\geq U_\tau).
\end{split}
\end{equation}
To see $(\dag)$ one can argue as follows: Suppose $\tau_p\leq n$ and
$S_{p,\tau_p}\geq U_{\tau_p}$.  Then $S_{q,\tau_p}\geq
S_{p,\tau_p}\geq U_{\tau_p}$. Hence, $\tau_q\leq \tau_p$.
Furthermore, for all $k\leq \tau_p$ we have $S_{q,k}\geq S_{p,k}>L_k$.
Hence, $S_{q,\tau_q}\geq U_{\tau_q}$.  Letting $n\to \infty$ in
(\ref{eq:ineqinlemma}) finishes the proof of the first inequality.
The second inequality can be shown similarly.
\end{proof}

\begin{proof}[Proof of Theorem \ref{th:stopinfinitetime}]
Suppose $L_n\geq U_n$ for some $n$. Then, by the definition of $U_n$ and $L_n$,
\begin{align*}
1=\Prob{\alpha}(\tau\leq n)
\leq
\Prob{\alpha}(S_\tau\geq U_\tau,\tau\leq n)  
+\Prob{\alpha}(S_\tau\leq L_\tau,\tau\leq n)
\leq 2\epsilon_n<2\epsilon<1,
\end{align*}
which is a contradiction. Hence, $U_n>L_n$.

Let $j_n=\left\lceil\Delta_n+n\alpha\right\rceil $,
where $\Delta_n=\sqrt{-n\log(\epsilon_n-\epsilon_{n-1})/2}$ and $\left\lceil x \right\rceil$ denotes the smallest integer greater than $x$.
We show $U_n\leq j_n$.
By a special case of Hoeffding's inequality,
 see \citet[Theorem 1]{okamoto1958sir} or \citet[Theorem 1]{hoeffding1963pis},
\begin{equation}
\label{eq:ineqintheorem}
\begin{split}
 \Prob{\alpha}(S_n\geq j_n, \tau\geq n)
&\leq \Prob{\alpha}(S_n\geq j_n)=\Prob{\alpha}(S_n/n-\alpha\geq j_n/n-\alpha)\\
&\leq\exp\left(-2n(j_n/n-\alpha)^2\right)\\
&\leq\exp\left(-2n\left(\frac{\Delta_n}{n}\right)^2\right)
=\epsilon_n-\epsilon_{n-1}.
\end{split}
\end{equation}
By the definition of $U_{n-1}$,
$$\Prob{\alpha}(S_\tau\geq U_\tau,\tau<n)
=\Prob{\alpha}(S_{n-1}\geq U_{n-1},\tau=n-1)
+\Prob{\alpha}(S_{\tau}\geq U_\tau,\tau<n-1)
 \leq \epsilon_{n-1}.$$
This together with the definition of $U_n$ and (\ref{eq:ineqintheorem}) yields $U_n\leq j_n$. Thus,
\begin{align}
\label{eq:asymptboundUn}
  \frac{U_n-n\alpha}{n}&\leq\frac{\Delta_n+1}{n}\to 0 \quad(n\to \infty)
\end{align}
Similarly, one can show 
\begin{equation}
  \label{eq:asymptboundLn}
  \frac{L_n-n\alpha}{n}\geq -\frac{\Delta_n+1}{n}\to 0 \quad(n\to \infty).
\end{equation}
Together with $U_n\geq L_n$ we get
 $U_n-n\alpha=o(n)$ and  $n\alpha-L_n=o(n)$.

Next, we show $\Expec{p}(\tau) <\infty $ for $p> \alpha$.
For large $n$ (say $n\geq n_0$),
$$
p-\frac{U_n}{n}=(p-\alpha)-(\frac{U_n}{n}-\alpha)=(p-\alpha)+o(1)\geq (p-\alpha)/2.
$$
For $n\geq n_0$, since $U_n/n-p\leq 0$,  Hoeffding's inequality shows  
  \begin{align*}
    \Prob{p}(\tau>n)\leq \Prob{p}(S_n<U_n)
=\Prob{p}\left(\frac{S_n}{n}-p\leq \frac{U_n}{n}-p\right)
\leq \exp\left(-2n\left(\frac{U_n}{n}-p\right)^2\right).
  \end{align*}
Thus $\Prob{p}(\tau>n)\leq \exp(-n(p-\alpha)^2/2)
$ for $n\geq n_0$.
Hence, 
\begin{align*}
  \Expec{p}(\tau)=&\int_0^\infty \Prob{p}(\tau>n)dn
\leq n_0+\int_{n_o}^\infty\exp(-n(p-\alpha)^2/2)dn\\
=&n_0+\frac{2}{(p-\alpha)^2}\exp(-n_0(p-\alpha)^2/2)
<\infty.
\end{align*}
Similarly, one can show  $\Expec{p}(\tau) <\infty $ for  $p< \alpha$.

To see (\ref{eq:boundhitwrongboundary}):
Let  $p\leq\alpha$. By the previous lemma we have
$\Prob{p}(\tau<\infty,S_\tau\geq U_\tau) \leq 
\Prob{\alpha}(\tau<\infty,S_\tau\geq U_\tau)$.
Since $\epsilon_n < \epsilon$, equation (\ref{eq:hitputuntilk}) implies $\Prob{\alpha}( \tau<\infty, S_\tau\geq U_\tau) \leq \epsilon$.  
Thus   $\Prob{p}(\tau<\infty,S_\tau\geq U_\tau) \leq \epsilon$.
The second part of  (\ref{eq:boundhitwrongboundary})  can be shown similarly.
\end{proof}

\begin{proof}[Proof of Theorem \ref{th:MLE}]
By (\ref{eq:boundhitwrongboundary})
it suffices to show that $\tau<\infty$, $S_\tau \geq U_\tau$ implies $\hat p>\alpha$ 
and that $\tau<\infty$, $S_\tau \leq L_\tau$ implies $\hat p\leq\alpha$.

First, we  show $L_n<n\alpha$.
Suppose $L_n\geq n\alpha$. 
Let $c=\lceil n\alpha \rceil$. 
Hence, by \cite[Satz 6]{Uhlmann1966Vdh}, if
$c\leq\frac{n-1}{2}$,
$$
\frac{1}{2}\leq \Prob{\frac{c}{n-1}}(S_n\leq c)\leq \Prob{\alpha}(S_n\leq c),
$$
since  $\frac{c}{n-1}\geq \frac{n\alpha}{n-1}=\alpha+\frac{\alpha}{n-1}>\alpha$;
and if $c\geq \frac{n-1}{2}$,
$$
\frac{1}{2}\leq \Prob{\frac{c+1}{n+1}}(S_n\leq c)\leq \Prob{\alpha}(S_n\leq c).
$$
since  $\frac{c+1}{n+1}\geq \frac{n\alpha+1}{n+1}=\alpha+\frac{1-\alpha}{n+1}>\alpha$.
Hence, 
\begin{align*}
  \frac{1}{2}&\leq \Prob{\alpha} (S_n\leq c)
\leq  \Prob{\alpha} (S_n\leq L_n)
\leq  \Prob{\alpha} (\tau\leq n)
\leq 2\epsilon_n <\frac{1}{2}
\end{align*}
which is a contradiction. 
$U_n\leq n\alpha$  leads to a contradiction in a similar way.
\end{proof}

\bibliographystyle{apalike}
\bibliography{lit,../../../lit/archive}

\end{document}